\documentclass[12pt]{article}
\usepackage{amsmath,amsfonts,amssymb,latexsym,amsthm}
\usepackage[a4paper,portrait]{geometry}
\usepackage{url}

\newcommand{\EM}{\ensuremath}

\newtheorem{thm}{Theorem}[section]%
\newtheorem{cor}[thm]{Corollary}%
\newtheorem{lem}[thm]{Lemma}%
\newtheorem{xpl}[thm]{Example}%
\newtheorem{rem}[thm]{Remark}%

\newcommand{\dN}{\EM{\mathbb{N}}}

\newcommand{\dP}{\EM{\mathbb{P}}}

\newcommand{\dR}{\EM{\mathbb{R}}}

\newcommand{\rL}{\EM{\mathrm{L}}}

\newcommand{\cC}{\EM{\mathcal{C}}}

\newcommand{\cF}{\EM{\mathcal{F}}}

\newcommand{\cL}{\EM{\mathcal{L}}}
\newcommand{\cM}{\EM{\mathcal{M}}}

\newcommand{\cX}{\EM{\mathcal{X}}}

\newcommand{\cZ}{\EM{\mathcal{Z}}}

\newcommand{\al}{\alpha}

\newcommand{\de}{\delta}

\newcommand{\la}{\lambda}

\newcommand{\Om}{\Omega}
\newcommand{\om}{\omega}

\newcommand{\si}{\sigma}
\newcommand{\Te}{\Theta}
\newcommand{\te}{\theta}

\newcommand{\p}[4]{{#3}\!\left#1{#4}\right#2} 
\newcommand{\NRM}[1]{\EM{{\left\| #1\right\|}}} 
\newcommand{\PAR}[1]{\EM{{\left(#1\right)}}} 
\newcommand{\entf}[1]{\mathbf{Ent}_{#1}}
\newcommand{\ent}[2]{\p(){\entf{#1}}{#2}}
\newcommand{\OL}[1]{\overline{#1}}
\newcommand{\WH}[1]{\widehat{#1}}
\newcommand{\HYP}[1]{\textbf{(A\ref{hyp:#1})}}

\title{\textsf{On the strong consistency \\ of asymptotic M-estimators}}
\author{Djalil~\textsc{Chafa\"\i} and Didier~\textsc{Concordet}}
\date{{\small\textsf{%
 Preprint. September 2006. \\
 Accepted for publication in Journal of Statistical Planning and Inference.}}}

\begin{document}

\maketitle

\begin{abstract}
  The aim of this article is to simplify Pfanzagl's proof of consistency for
  asymptotic maximum likelihood estimators, and to extend it to more general 
  asymptotic $M$-estimators. 
  The method relies on the existence of a sort of contraction
  of the parameter space which admits the true parameter as a fixed point. The
  proofs are short and elementary.
\end{abstract}

\section{Introduction}

After the seminal work\footnote{The interested reader may find a quite recent
  account in \cite{MR1617519} and references therein.} of Fisher, the
asymptotic properties of maximum likelihood estimators, and in particular
their consistency, were studied by various authors, including Doob
\cite{MR1501765}, Cram\'er \cite{MR0016588}, and Huzurbazar \cite{MR0028000}.
Nowadays, one of the best known result regarding consistency goes back to
Wald, who gave in \cite{MR0032169} a short and elegant proof of strong
consistency of parametric maximum likelihood estimators. Since that time,
several authors studied various versions of such consistency problems,
including among others, Le Cam \cite{MR0054913}, Kiefer and Wolfowitz
\cite{MR0086464}, Bahadur \cite{MR0207085,MR0315820}, Huber \cite{MR0216620},
Perlman \cite{MR0405668}, Wang \cite{MR803749}, and Pfanzagl
\cite{MR944202,MR1074441}.

Wald's original proof relies roughly on local compactness of the parameter
space, on continuity and coercivity\footnote{By coercivity we mean that the
  log-likelihood tends to $-\infty$ when the parameter tends to $\infty$.} of
the log-likelihood, on the law of large numbers, and last but not least on
local uniform integrability of the log-likelihood. It does not require
differentiability, and makes extensive use of likelihood ratios. The
integrability assumption has been weakened by many authors, including for
instance Kiefer and Wolfowitz in \cite{MR0086464} and Perlman in
\cite{MR0405668}, see also \cite{MR0315820}. One can find a modern
presentation of Wald's method for $M$-estimators in van der Vaart's monograph
\cite{MR1652247}.

Pfanzagl gave in \cite{MR944202,MR1074441} a proof of strong consistency of
asymptotic maximum likelihood estimators for nonparametric ``concave models''
with respect to the estimated parameter, including nonparametric mixtures. His
approach relies in particular on a simplification of an earlier work of Wang
in \cite{MR803749} based on uniform local bound of the likelihood ratio. 

The present work was initially motivated by the inverse problems
considered in \cite{chafai-loubes}. Our aim is to simplify Pfanzagl's 
approach, and to extend the framework from asymptotic maximum likelihood to 
more general asymptotic M-estimators. 
In particular, log-likelihood ratios are replaced by contrast
differences. The hypotheses appearing in our main Theorem are unnecessarily
strong. However, they allow a simple and short presentation. We emphasize the
role played by a sort of contraction map $a^*$ defined on the parameter space.
We do not assume any coercivity of the contrast as in \cite{MR0032169}.
However, we require the compactness of the space of the estimated parameter,
as in \cite{MR0086464} and \cite{MR1652247} for example. This compactness
comes usually for free in the case of fully nonparametric models. We do not
make use of any Uniform Law of Large Numbers. Our method does not belong to
the Glivenko-Cantelli approaches of consistency, as in \cite{MR1652319},
\cite{MR1797759}, \cite{MR1365666}, \cite{MR1652247} and
\cite{MR1973724,MR1739079} and references therein.

Let $\Te$ be a separable Hausdorff topological space with countable base. Let
$(P_\te)_{\te\in\Te}$ be a known family of Borel measures on a measurable
space $\cX$. Let $\te^*\in\Te$ be some unknown point of $\Te$ such that
$P^*:=P_{\te^*}$ is a probability measure. Let $(X_n)_{n\in\dN}$ be an i.i.d.
sequence of observed random variables defined on a probability space
$(\Om,\cF,\dP)$ and taking their values in $\cX$, with common law $P^*$. Let
$(\WH{\te}_n)_{n\in\dN}$ be a sequence of random variables defined on
$(\Om,\cF,\dP)$, taking their values in $\Te$, and such that
$(\WH{\te}_n)_{n\in\dN}$ is $\cF_n$-measurable for any $n\in\dN$, where
$\cF_n:=\si(X_0,\ldots,X_n)$. We say that $(\WH{\te}_n)_{n\in\dN}$ is
\emph{strongly consistent} if and only if
\begin{equation}\label{eq:strong-consis}
  \dP-\text{a.s.\quad }\lim_{n\to+\infty}\WH{\te}_n = \te^*.
\end{equation}
We use in the sequel the abbreviations ``a.s.'' for \emph{almost sure},
``a.a.'' for \emph{almost all}, and ``a.e.'' for \emph{almost everywhere}. Let
$\Te\times\cX\ni(\te,x)\mapsto m(\te,x)\in\dR$ be a known function such that
$m_\te:=m(\te,\cdot)$ is measurable for any $\te\in\Te$. For any $n$, we
define the random function $M_n:\Te\to\dR$ by
\begin{equation*}
  M_n(\te):=\frac{1}{n}\sum_{i=1}^n m(\te,X_i).
\end{equation*}
This can be written also $M_n(\te)=\dP_n m_\te$ where
$\dP_n:=\frac{1}{n}(\de_{X_1}+\cdots+\de_{X_n})$ is the empirical measure. We
say that $(\WH{\te}_n)_n$ is a \emph{sequence of asymptotic M-estimators} if
and only if
\begin{equation}\label{eq:ame}
  \dP-\text{a.s.\ }
  \varlimsup_{n\to+\infty}\PAR{\sup_{\Te}M_n-M_n(\WH{\te}_n)}=0.
\end{equation}
The term \emph{asymptotic} is used for the same notion (with the likelihood)
by Pfanzagl in \cite{MR944202}. In the literature, some authors, including
Wald and Perlman, use the term \emph{approximate} rather than
\emph{asymptotic}. However, the term \emph{approximate} has been used by
Bahadur in a different sense in \cite[page 34]{MR0315820}.

For example, if for large enough $n$, there exists an $\cF_n$-measurable
$\WH{\te}_n$ in $\Te$ such that $M_n(\WH{\te_n})=\sup_\Te M_n$, then such a
random sequence $(\WH{\te}_n)_{n\in\dN}$ fulfils \eqref{eq:ame}.

For any probability measure $P$ on $\cX$, let $\rL_+^1(\cX,P)$ (resp.
$\rL_-^1(\cX,P)$) be the set of random variables $Z:\cX\to\dR$ such that
$Z^+:=\max(+Z,0)$ (resp. $Z^-:=\max(-Z,0)$) is in $\rL^1(\cX,P)$. On
$E(\cX,P):=\rL_-^1(\cX,P)\cup \rL_+^1(\cX,P)$, the expectation
$P(Z)=P(Z^+)-P(Z^-)$ makes sense and takes its values in
$\OL{\dR}:=\dR\cup\{\pm\infty\}$. For any $\te\in\Te$ such that $m_\te\in
E(\cX,P^*)$, we define the \emph{contrast} $M^*(\te)\in\OL{\dR}$ by
\begin{equation}\label{eq:def-M-star}
  M^*(\te):=P^* m_\te.
\end{equation}
In the sequel, we say that the model is \emph{identifiable} when for any
$\te\in\Te$, the condition $P_\te=P^*$ implies that $\te=\te^*$.

\begin{xpl}[Log-Likelihood]
  Assume that for some fixed Borel measure $Q$ on $\cX$, one has $P_\te\ll Q$
  for any $\te\in\Te$. Let $f_\te:=dP_\te/dQ$ and assume that $f_\te>0$ on
  $\cX$ for any $\te\in\Te$. Define $m(\te,x):=\log(f_\te(x))$. Then
  $M_n:\Te\to\dR$ is the \emph{log-likelihood} random functional given by
  $M_n(\te)=\dP_n m_\te=\dP_n\log(f_\te)$. We will speak about sequences of
  ``asymptotic maximum likelihood estimators''. The log-likelihood ratio is
  $$
  M_n(\te_1)-M_n(\te_2)=\dP_n\log(f_{\te_1}/f_{\te_2}).
  $$
  As usual for the log-likelihood, when $M^*(\te^*)$ is finite, one can write
  for any $\te$
  $$
  M^*(\te)-M^*(\te^*)=-\ent{}{P_{\te^*}\,\vert\,P_{\te}},
  $$ 
  where $\ent{}{P_{\te_1}\,\vert\,P_{\te_2}}$ is the Kullback-Leibler relative
  entropy of $P_{\te_1}$ with respect to $P_{\te_2}$. In particular,
  $M^*(\te)\leq M^*(\te^*)$ with equality if and only if $P_{\te}=P_{\te^*}$,
  which implies $\te=\te^*$ if the model is identifiable. Notice that when $Q$
  is the Lebesgue measure on $\cX=\dR^n$, then $-M^*(\te^*)=-\int_\cX
  f_{\te^*}(x)\log(f_{\te^*}(x))\,dx$ is the Shannon entropy of $f_{\te^*}$.
\end{xpl}

\begin{xpl}[Beyond the log-likelihood]
  Assume that for some fixed Borel measure $Q$ on $\cX$, one has $P_\te\ll Q$
  for any $\te\in\Te$, with $P_\te(\cX)\leq1$ and $f_\te:=dP_\te/dQ$. Let
  $\Phi,\Psi:(0,+\infty)\to\dR$ be two smooth functions. Assume that
  $\Psi(f_\te)\in\rL^1(\cX,Q)$ for any $\te\in\Te$. Define $m_\te$ by
  $$
  m_\te = \Phi(f_\te) - \int_\cX\!\Psi(f_\te)\,dQ + P_\te(\cX).
  $$
  This gives rise the the following empirical contrast
  $$
  M_n(\te) = \dP_n(\Phi(f_\te))-\int_\cX\!\Psi(f_\te)\,dQ +P_\te(\cX).
  $$
  In particular, if $\te\in\Te$ is such that $\Phi(f_\te)\in\rL^1(\cX,P^*)$
  where here again $P^*:=P_{\te^*}$,
  $$
  M^*(\te) = P^*(\Phi(f_\te))-\int_\cX\!\Psi(f_\te)\,dQ + P_\te(\cX).
  $$
  Assume now that $u\mapsto u\Phi'(u)$ is locally integrable on $\dR_+$, and
  consider the case where $\Psi$ is the $\Phi$-transform given for any
  $u\in(0,+\infty)$ by
  $$
  \Psi(u)=\int_0^u\!v\Phi'(v)\,dv.
  $$
  For $\Phi:u\mapsto\log(u)$, one has $\Psi:u\mapsto u$ and we recover the
  log-likelihood contrast
  $$
  M^*(\te)=P^*(\log(f_\te)).
  $$ 
  For $\Phi:u\mapsto u$, one has $\Psi:u\mapsto \frac{1}{2}u^2$, and we get
  the quadratic contrast
  $$
  M^*(\te) %
  =-\mbox{$\frac{1}{2}$}\NRM{f_\te-f_{\te^*}}^2_{\rL^2(\cX,Q)} %
  +\mbox{$\frac{1}{2}$}\NRM{f_{\te^*}}^2_{\rL^2(\cX,Q)} %
  +P_\te(\cX).
  $$
  In both cases, the map $\te\mapsto M^*(\te)$ admits $\te^*$ as unique
  maximum provided that the model is identifiable. More generally, 
  define the $\Phi$-transform $\Te:(0,+\infty)^2\to\dR$ by
  \begin{align*}
  \Te(u,v) :&=u\Phi(v)-\Psi(v) \\
  &=u\Phi(v)-\int_0^v\!w\Phi'(w)\,dw. 
  \end{align*}
  When $\te$ and $\te^*$ are such that both $\Te(f_{\te^*},f_{\te^*})$ and
  $\Te(f_{\te^*},f_\te)$ belong to $\rL^1(\cX,Q)$,
  $$
  M^*(\te) = \int_\cX\!(\Te(f_{\te^*},f_\te)-\Te(f_{\te^*},f_{\te^*}))\,dQ %
  + \int_\cX\!\Te(f_{\te^*},f_{\te^*})\,dQ %
  + P_\te(\cX).
  $$
  Notice that $\Te$ is linear in $\Phi$. One can consider useful examples for
  which the function $\Phi$ is bounded, in such a way that $m_\te$ is bounded
  for any $\te\in\Te$. For instance, let us examine the case where
  $\Phi:u\mapsto-(1+u)^{-2}$. Then, $\Psi:u\mapsto -u^2(1+u)^{-2}$, and the
  map $\te\mapsto M^*(\te)$ admits $\te^*$ as unique maximum, provided
  identifiability holds, since for any $(u,v)\in\dR_+^2$,
  $$
  \Te(u,v)=-\frac{u+v^2}{(1+v)^2} %
  \text{\quad and \quad} %
  \Te(u,v)-\Te(u,u)=-\frac{(v-u)^2}{(1+u)(1+v)^2}.
  $$
  The function $\Psi$ is additionally bounded here. The similar case
  $\Phi:u\mapsto-(1+u^2)^{-1}$ is also quite interesting. Notice that
  $\Te(u,\cdot)$ is concave on $(0,+\infty)$ as soon as $\Phi$ is concave, non
  decreasing, with $\Phi'(v)+v\Phi''(v)\geq0$ for any $v>0$. Observe that this
  is not the approach of Pfanzagl in \cite{MR1074441}, which is more related
  to the log-likelihood ratio. Notice that in the case of the log-likelihood,
  one has $\Phi:u\mapsto\log(u)$, which gives $\Psi:u\mapsto u$ and
  $\Te:(u,v)\mapsto-u\log(v)-v$, and thus $\Te(u,v)-\Te(u,u)=u\log(u/v)+u-v$.
  It might be possible to extensively study such ``$\Phi$-estimators'', in the
  spirit of the ``$\Phi$-calculus'' developed in \cite{MR2081075,chafai-mmi}.
  This is however outside the scope of this short article.
\end{xpl}

One can notice that the observation of Lindsay in \cite{MR684866,MR707929}
regarding the nature of maximum likelihood for nonparametric mixture models
remains valid for more general models provided that $m$ is concave.

\section{Main result and Corollaries}

With the settings given in the Introduction, the following Theorem holds.

\begin{thm}\label{th:consistency-a-la-pf}
  Assume that $\Te$ is compact and that the following assumptions hold.
  \begin{enumerate}
    \renewcommand{\labelenumi}{\textbf{\quad(A\theenumi)}}
  \item\label{hyp:cont} For $P^*$-a.a. $x\in\cX$, the map $m(\cdot,x)$ is
    continuous on $\Te$;
  \item\label{hyp:sep} There exists a continuous map $a^*:\Te\to\Te$ which may
    depend on $\te^*$ such that for any $\te\neq\te^*$, there exists a
    neighborhood $V\subset\Te$ of $\te$ for which $\sup_V\PAR{m-m_{a^*}}\in
    \rL_+^1(\cX,P^*)$ and $P^*(m_\te-m_{a^*(\te)})<0$.
  \end{enumerate}
  Then any sequence $(\WH{\te}_n)_n$ of asymptotic M-estimators is strongly
  consistent.
\end{thm}

\begin{proof}
  Postponed to section \ref{se:proof:th:consistency-a-la-pf}.
\end{proof}

The quantity $P^*(m_\te-m_{a^*(\te)})$ in \HYP{sep} has a meaning in
$\OL{\dR}$ since the first part of \HYP{sep} ensures that
$m_\te-m_{a^*(\te)}\in\rL_+^1(\cX,P^*)$. Moreover, $P^*(m_\te-m_{a^*(\te)})$
reads $M^*(\te)-M^*(a^*(\te))$ when the couple $(m_\te,m_{a^*(\te)})$ is in
$\rL^1_-(\cX,P^*)\times\rL^1_+(\cX,P^*)$ or in
$\rL^1_+(\cX,P^*)\times\rL^1_-(\cX,P^*)$.

Since $\te^*$ is unknown in practice, each assumption in Theorem
\ref{th:consistency-a-la-pf} must hold for any $\te^*\in\Te$ such that
$P_{\te^*}$ is a probability measure, in order to make the result useful.

\begin{rem}[Assumptions]
  The first part of \HYP{sep} is in a way an $M$-estimator version of the
  integrability condition considered by Kiefer and Wolfowitz for the
  log-likelihood in \cite{MR0086464}. The assumptions \HYP{cont} and \HYP{sep}
  required by Theorem \ref{th:consistency-a-la-pf} can be weakened. However,
  they permit a streamlined presentation. In particular, only lower
  semi-continuity is needed in \HYP{cont}, see for instance \cite{MR944202}.
  Additionally, and following for example \cite[page 266]{MR0405668}, the
  uniform integrability assumption \HYP{sep} can be weakened, by considering
  blocks of $k>1$ observations instead of one observation, see also
  \cite[comments following Theorem 5.14]{MR1652247}.
\end{rem}

As stated in the following Corollary, Theorem \ref{th:consistency-a-la-pf}
implies a version of Wald consistency Theorem for asymptotic $M$-estimators,
see \cite{MR0032169}, \cite[Section 2 page 269]{MR0405668}, and \cite[Theorem
5.14]{MR1652247}.

\begin{cor}[Perlman-Wald]\label{co:wald}
  Assume that $\Te$ is compact, and that for $P^*$-a.a. $x\in\cX$, the map
  $m(\cdot,x)$ is continuous on $\Te$. Assume that for any $\te$ in $\Te$,
  there exists a neighborhood $V$ such that $\sup_V m\in\rL^1(\cX,P^*)$.
  Assume in addition that $M^*$ achieves its supremum over $\Te$ at $\te^*$,
  and only at $\te^*$. Then, any sequence of asymptotic $M$-estimators is
  strongly consistent.
\end{cor}

\begin{proof}
  One has $m_\te\in\rL^1(\cX,P^*)$ for any $\te$ in $\Te$, and thus
  $M^*:\Te\to\dR$ is well defined. Moreover, \HYP{sep} holds with a constant
  map $a^*\equiv\te^*$. Namely, for any $\te\neq\te^*$, one has on one hand
  $P^*(m_\te-m_{\te^*})<0$ since $M^*(\te)<M^*(\te^*)$, and on the other hand
  $$
  \sup_V(m-m_{a^*})=-m_{\te^*}+\sup_V m\in\rL^1(\cX,P^*).
  $$
\end{proof}

As stated in the following Corollary, Theorem \ref{th:consistency-a-la-pf}
implies the main result of Pfanzagl in \cite{MR944202} for concave models,
itself based on an earlier result of Wang in \cite{MR803749}. This is
typically the case for mixtures models, for which $\Te$ is a convex set of
probability measures on some measurable space, cf. section \ref{se:mix}.

\begin{cor}[Pfanzagl-Wang]\label{co:pf}
  Let $Q$ be a reference Borel measure on $\cX$. Consider the case where $\Te$
  is a convex compact subset of a linear space such that for any $\te\in\Te$,
  $P_\te(\cX)\leq1$ and $P_\te\ll Q$ with $f_\te:=dP_\te/dQ>0$ on $\cX$.
  Suppose that $Q$-a.e. on $\cX$, the map $\te\mapsto f_\te(x)$ is concave and
  continuous on $\Te$. Assume that the model is identifiable. Consider
  $m_\te:=\log(f_\te)$ and the related log-likelihood $M_n$. Then any sequence
  of asymptotic log-likelihood estimators is strongly consistent.
\end{cor}

\begin{proof}
  First of all, we notice that it is not possible to take $a^*\equiv\te^*$
  since we cannot ensure that the condition
  $m_{\te^*}-m_\te=\log(f_{\te^*}/f_\te)\in\rL_+^1(\cX,P^*)$ of \HYP{sep} is
  true. However, the concavity of the model allows to take a map $a^*$ which
  is a strict contraction around $\te^*$. Namely, for an arbitrary
  $\la\in(0,1)$, let us take
  $$
  a^*(\te):=\la\te^*+(1-\la)\te.
  $$
  The concavity of the model yields
  $$
  m_{a^*(\te)}-m_\te %
  =\log\PAR{\frac{f_{\la\te^*+(1-\la)\te}}{f_\te}} %
  \geq\log\PAR{\frac{\la f_{\te^*}+(1-\la)f_\te}{f_\te}} %
  \geq \log(1-\la).
  $$
  Now, we have $\log(1-\la)\in\rL^1(\cX,P^*)$ since $\la<1$. Define the
  function $\Phi:\dR_+\to\dR$ by $\Phi(u):=u\log(\la u + (1-\la))$. The
  concavity of the model yields
  $$
  P^*(m_{a^*(\te)}-m_\te) %
  \geq\int_\cX f_{\te^*} %
  \log\PAR{\frac{\la f_{\te^*}+(1-\la)f_\te}{f_\te}}\,dQ %
  =\int_\cX \Phi\PAR{\frac{f_{\te^*}}{f_\te}}f_\te\,dQ.
  $$
  Let us show that the right hand side of the inequality above is strictly
  positive when $\te\neq\te^*$. One has $P_\te(\cX)>0$ since $f_\te>0$. Define
  $\Psi(u):=u\Phi(1/u)$. Jensen's inequality for the probability measure
  $P_\te(\cX)^{-1}P_\te$ and the convex function $\Phi$ yields
  \begin{equation}\label{eq:J}
    \int_\cX \Phi\PAR{\frac{f_{\te^*}}{f_\te}}f_\te\,dQ %
    \geq \Psi(P_\te(\cX)).
  \end{equation}
  It is enough to show that either \eqref{eq:J} is strict or the right hand
  side of \eqref{eq:J} is strictly positive. Since $\la>0$, the function
  $\Phi$ is strictly convex. Thus equality holds in \eqref{eq:J} if and only
  if $P_\te(f_{\te^*}=\al f_\te)=1$ for some $\al\in\dR_+$. The only
  admissible case is $\al=P_\te(\cX)^{-1}>1$ since $P_{\te^*}(\cX)=1$ and
  since identifiability forbids $P_\te(f_{\te^*}=f_\te)=1$. Therefore, if
  $P_\te(\cX)=1$, inequality \eqref{eq:J} is necessarily strict. On the other
  hand, $\Psi(1)=0$ and $\Psi(u)>0$ when $u<1$. Thus the right hand side of
  \eqref{eq:J} is always non negative, and is strictly positive as soon as
  $P_\te(\cX)<1$. We conclude that $P^*(m_{a^*(\te)}-m_\te)>0$ as soon as
  $\te\neq\te^*$. This shows that \HYP{sep} holds with $V=\Te$, and the proof
  is thus complete.
\end{proof}

\begin{rem}[About the map $a^*$]\label{rm:a-xtrm}
  Let $a^*:\Te\to\Te$ be a map which satisfies the condition
  $P^*(m_\te-m_{a*(\te)})<0$ for any $\te\neq\te^*$ of \HYP{sep}. Then, the
  impossibility of $P^*(m_\te-m_\te)<0$ for any $\te$ yields that
  \begin{itemize}
  \item $a^*(\te)\neq\te$ for any $\te\neq\te^*$. In particular,
    \begin{itemize}
    \item the map $a^*$ cannot be the identity map ;
    \item if $a^*$ is constant, then $a^*\equiv\te^*$ ;
    \item the point $\te^*$ is the only possible fixed point for $a^*$.
    \end{itemize}
  \end{itemize}
  The proof of Corollary \ref{co:wald} gives an example where $a^*\equiv\te^*$
  works and fulfills \HYP{sep}. In contrast, Corollary \ref{co:pf} provides a
  situation where a constant $a^*$ does not fulfill \HYP{sep}. However, we
  have shown in the proof of Corollary \ref{co:pf} that an $a^*$ map which is
  a strict contraction around $\te^*$ fulfills \HYP{sep}. Actually, when $\Te$
  has the structure of a convex subset of a vector space, any strict
  contraction around $\te^*$ fulfills the properties of $a^*$ listed above.
  The existence of a fixed point can be related to Brouwer-like fixed point
  Theorems. For instance, any continuous mapping of a non-empty compact convex
  subset of $\dR^d$ into itself contains at least one fixed point.
  Consequently, when $\Te$ is a non-empty compact and convex subset of
  $\dR^d$, any continuous $a^*$ map admits $\te^*$ as a unique fixed point.
  There exists numerous dimension free Brouwer-like fixed points theorems, due
  to Schauder, Tikhonov, Kakutani, \ldots, see for instance \cite{MR816732}
  and \cite{MR1996163}.
\end{rem}

\begin{rem}[Infinite values of $m$]\label{rm:m-inf}
  Theorem \ref{th:consistency-a-la-pf} does not allow $m$ to take the value
  $-\infty$. This limitation is due to the fact that differences of the form
  $m_\te-m_{\te'}$ do not make sense if $m$ is allowed to take the value
  $-\infty$. The consistency proof of Wald does not suffer from such a
  limitation since it does not rely on $m$ differences, but it requires
  however strong uniform integrability assumptions. A careful reading of the
  proof of Theorem \ref{th:consistency-a-la-pf} shows that only differences of
  the form $m_\te-m_{a^*(\te)}$ are involved. On the other hand, according to
  Remark \ref{rm:a-xtrm}, $a^*(\te)\neq\te$ for any $\te\neq\te^*$.
  Consequently, one may allow, in Theorem \ref{th:consistency-a-la-pf}, the
  map $m(\te,x)$ to take the value $-\infty$ for at most one value of $\te$.
  For the log-likelihood, $m_\te=\log(f_\te)$ and one has $m_\te(x)=-\infty$
  if and only if $f_\te(x)=0$. One may allow $f_\te\equiv0$ for at most one
  value of $\te$ in Corollary \ref{co:pf}.
\end{rem}

\begin{rem}
  Let $\te\in\Te$ such that $m_\te\in E(\cX,P^*)$. Then, the law of large
  numbers applies and gives that $P^*$-a.s., $\lim_n
  M_n(\te)=M^*(\te)\in\OL{\dR}$, and the a.s. subset of $\cX$ may depend on
  $\te$. In particular $M_n(\te)=M^*(\te)+o_P(1)$. For a sequence
  $(\WH{\te}_n)_n$ satisfying \eqref{eq:ame}, one can write for any
  $\te\in\Te$ with finite $M_n(\te)$
  \begin{align*}
    M_n(\WH{\te}_n)
    &= M_n(\WH{\te}_n)-M_n(\te) + M_n(\te) \\
    &\geq -\PAR{\sup_\Te M_n-M_n(\WH{\te}_n)} + M_n(\te) \\
    &= o_P(1)+ M(\te)
  \end{align*}
  where the last step follows by \eqref{eq:ame} and the law of large numbers.
\end{rem}

\section{Log-Likelihood and mixtures models}
\label{se:mix}
 
For any topological space $\cZ$ equipped with its Borel $\si$-field, we denote
by $\cM_1(\cZ)$ the set of probability measures on $\cZ$, and by $\cC_b(\cZ)$
the set of bounded real valued continuous functions on $\cZ$. The Prohorov
topology on $\cM_1(\cZ)$ is defined as follows: $\te_n\to\te$ in $\cM_1(\cZ)$
if and only if $\int_\cZ\!f\,d\te_n\to\int_\cZ\!f\,d\te$ for any
$f\in\cC_b(\cZ)$. It is known that a subset of $\cM_1(\cZ)$ is compact if and
only if it is tight. As a consequence, $\cM_1(\cZ)$ is not compact in general.
Following \cite[section 5 page 149]{MR944202}, the set sub-probabilities
provides a compactification which allows the following consistency result for
asymptotic log-likelihood estimators of nonparametric mixture models.

\begin{cor}[Pfanzagl]\label{co:pfm}
  Let $\cZ$ be a locally compact Hausdorff topological space with countable
  base. Let $Q$ be a measure on a measurable space $\cX$. Let
  $k:\cX\times\cZ\to(0,+\infty)$ be such that $\int\!k(x,z)dQ(x)=1$ for any
  $z\in\cZ$ and $k(x,\cdot)\in\cC_b(\cZ)$ for any $x\in\cX$. Let
  $\Te:=\cM_1(\cZ)$ and consider the family $(P_\te)_{\te\in\Te}$ of
  probability measures on $\cX$ defined by $dP_\te=f_\te dQ$ with
  $f_\te(x):=\int\!k(x,z)\,d\te(z)$. Assume that the model is identifiable.
  Let $m:\Te\times\cX\to\dR$ be the map defined by $m(\te,x):=\log f_\te(x)$,
  and $M_n$ be the corresponding log-likelihood. Then any sequence of
  asymptotic maximum likelihood estimators is strongly consistent for the
  Prohorov topology.
\end{cor}

\begin{proof}
  As explained above, $\Te=\cM_1(\cZ)$ is not compact for the Prohorov
  topology, and one must consider a suitable compactification, as in
  \cite{MR0315820} for instance. Let $\cC_0(\cZ)$ be the set of real valued
  continuous functions on $\cZ$ which vanish at infinity. Let $\OL{\Te}$ be
  the set of Borel measures $\te$ on $\cZ$ such that $\te(\cZ)\leq 1$ (i.e.
  sub-probabilities), equipped with the vague topology related to
  $\cC_0(\cZ)$. Namely, $\te_n\to\te$ in $\OL{\Te}$ if and only if
  $\int_\cZ\!f\,d\te_n\to\int_\cZ\!f\,d\te$ for any $f\in\cC_0(\cZ)$. The
  injection $\Te\subset\OL{\Te}$ is continuous; $\OL{\Te}$ is a compact
  metrizable topological space, and thus has a countable base. Moreover,
  $\OL{\Te}$ is convex, and for any $\te\in\OL{\Te}$, there exists
  $\te'\in\Te$ and $\al\in[0,1]$ such that $\te=\al\te'$.

  We extend the set of probability measures $(P_\te)_{\te\in\Te}$ on $\cX$ to
  the set of sub-probability measures $(P_\te)_{\te\in\OL{\Te}}$ on $\cX$,
  where $dP_\te=f_\te dQ$ and $f_\te(x):=\int\!k(x,z)\,d\te(z)$. One has by
  virtue of Fubini-Tonelli Theorem that $P_\te(\cX)=\te(\cZ)$, and thus
  $P_\te\in\cM_1(\cX)$ if and only if $\te\in\Te:=\cM_1(\cZ)$. Notice that
  $\te^*$ is taken in $\Te$.

  Let $\te\in\OL{\Te}$ such that $P_\te=P_{\te^*}$. Since $\te^*$ is taken in
  $\Te$, one has that $P_\te\in\cM_1(\cX)$, therefore $\te\in\Te$ and thus
  $\te=\te^*$ by identifiability in $\Te$. Notice that $\OL{\Te}$ is the
  convex envelope of $\Te\cup\{0\}$. The set $\OL{\Te}$ contains the null
  measure $0$, for which $f_0\equiv0$ and thus $m_0\equiv-\infty$. If
  $\te\in\OL{\Te}$ with $\te\neq0$, then $f_\te>0$ on $\cX$ since $k>0$, and
  thus $m_\te(x):=\log f_\te(x)$ is finite for any $x\in\cX$. For any
  $x\in\cX$, the map $\te\in\OL{\Te}\mapsto m_\te(x)$ is continuous since
  $k(x,\cdot)$ is in $\cC_0(\cZ)$.

  For any $\te\in\OL{\Te}$ with $\te\neq0$, one can write $\te=\al\te'$ with
  $\te'\in\Te$ and $\al:=\te(\cZ)\in[0,1]$. One has then $f_{\te}=\al
  f_{\te'}$ and thus $m_{\te}=\log\al+m_{\te'}$. Therefore,
  $$
  M_n(\te)=\log\al+M_n(\te')\leq M_n(\te').
  $$
  As a consequence, $\sup_{\te\in\Te}M_n(\te)=\sup_{\te\in\OL{\Te}}M_n(\te)$,
  and one may substitute $\Te$ by $\OL{\Te}$ in the definition \eqref{eq:ame}.
  Now, let $(\WH{\te}_n)_{n\in\dN}$ be a sequence in $\Te$ of asymptotic
  maximum likelihood estimators. Corollary \ref{co:pf} and Remark
  \ref{rm:m-inf} for $(P_\te)_{\te\in\OL{\Te}}$ apply and give the $P^*$-a.s.
  convergence for the vague topology of $(\WH{\te}_n)_{n\in\dN}$ towards
  $\te^*$. Since both the sequence and the limit are in $\Te$, the convergence
  holds for the Prohorov topology, and the desired result is established.
\end{proof}

\begin{rem}
  A mixture model can always be seen as a conditional model. The observed
  random variables $X$ with values in $\cX$ is the first component of the
  couple $(X,Z)$ with values in $\cX\times\cZ$. The component $Z$ is not
  observed. However, the conditional law $\cL(X\,\vert\,Z=z)$ is known, and
  has density $k(\cdot,z)$ with respect to $Q$ on $\cX$. If $\te=\cL(Z)$, then
  $\cL(X)$ has density $f_\te$ with respect to $Q$ on $\cX$.
\end{rem}

\section{Proof of main result}
\label{se:proof:th:consistency-a-la-pf}

\begin{lem}[Reformulation]
  The random sequence $(\WH{\te}_n)_n$ is a sequence of asymptotic
  M-estimators if and only if
  \begin{equation}\label{eq:ame-cns}  
    \dP\text{--a.s.},\quad \forall (\te_n)_n\in\Te^\dN,\quad
    \varlimsup_{n\to+\infty} \PAR{M_n(\te_n)-M_n(\WH{\te}_n)}\leq0.
  \end{equation}
\end{lem}

\begin{proof}
  The proof is done ``$\om$ by $\om$'', and the a.s. sets in \eqref{eq:ame}
  and \eqref{eq:ame-cns} are the same. Recall that $(\WH{\te}_n)_n$ is a
  sequence of asymptotic M-estimators if and only if \eqref{eq:ame} holds.
  Actually, the definition of the supremum gives
  $\sup_{\te\in\Te}M_n(\te)-M_n(\WH{\te}_n)\geq0$. Therefore, \eqref{eq:ame}
  is equivalent to
  \begin{equation}\label{eq:ame-alt}
    \dP\text{--a.s.},\quad %
    \varlimsup_{n\to+\infty} %
    \PAR{\sup_{\te\in\Te}M_n(\te)-M_n(\WH{\te}_n)}\leq 0.
  \end{equation}
  The Lemma is thus reduced to the equivalence between \eqref{eq:ame-alt} and
  \eqref{eq:ame-cns}. We begin by the proof of the implication
  \eqref{eq:ame-alt} $\Rightarrow$ \eqref{eq:ame-cns}. Let $A$ be some
  $\dP$--a.s. set such that \eqref{eq:ame-alt} holds. We proceed by fixing
  $\om\in A$. We hide the dependency on $\om$ in the notation of $M_n$ and
  $\WH{\te}_n$ to lightweight the expressions. Let $(\te_n)_n$ be a sequence
  in $\Te$. By definition of the supremum, we have
  $M_n(\te_n)\leq\sup_{\te\in\Te}M_n(\te)$. Thus, we get
  $$
  M_n(\te_n)-M_n(\WH{\te}_n) \leq \sup_{\te\in\Te}M_n(\te)-M_n(\WH{\te}_n).
  $$
  Taking the $\varlimsup_{n\to+\infty}$ of both sides and using
  \eqref{eq:ame-alt} provides the expected result \eqref{eq:ame-cns}. It
  remains to establish the implication \eqref{eq:ame-cns} $\Rightarrow$
  \eqref{eq:ame-alt}. Let $A$ be some $\dP$--a.s. set such that
  \eqref{eq:ame-cns} holds. Here again, we proceed by fixing $\om\in A$, and
  we hide the dependency on $\om$ in the notation of the random objects like
  $M_n$ and $\WH{\te}_n$. By definition of the supremum, there exists, for any
  $n$, an element $\te_n\in\Te$ such that
  $$
  \sup_{\te\in\Te} M_n(\te)-M_n(\te_n)-\frac{1}{n} \leq 0.
  $$
  Notice that $\te_n$ depends on $\om$ since $M_n$ depends on $\om$. This
  yields
  \begin{equation}\label{eq:ame-cns-aux}
    \varlimsup_{n\to+\infty}\PAR{\sup_{\te\in\Te} M_n(\te)-M_n(\te_n)} \leq 0.
  \end{equation}
  Now we write the telescopic sum
  $$
  \sup_{\te\in\Te}M_n(\te)-M_n(\WH{\te}_n)
  = \sup_{\te\in\Te}M_n(\te) -M_n(\te_n) +M_n(\te_n) -M_n(\WH{\te}_n),
  $$
  which gives
  \begin{multline*}
  \varlimsup_{n\to+\infty}\PAR{\sup_{\te\in\Te}M_n(\te)-M_n(\WH{\te}_n)} \\
  \leq \varlimsup_{n\to+\infty}\PAR{\sup_{\te\in\Te}M_n(\te) -M_n(\te_n)}
  +\varlimsup_{n\to+\infty}\PAR{M_n(\te_n) -M_n(\WH{\te}_n)}. 
  \end{multline*}
  The two terms of the right hand side are ``$\leq 0$'' by virtue of
  \eqref{eq:ame-cns-aux} and \eqref{eq:ame-cns} respectively. This provides
  the desired result \eqref{eq:ame-alt}, as expected.
\end{proof}  

\begin{lem}[Separation]\label{le:consistency-from-sep-seq}
  Assume that $\dP$--a.s., for any neighborhood $U$ of $\te^*$, for any
  sequence $(\te_n)_n$ in $U^c$, there exists a sequence $(\te'_n)_n$ in $\Te$
  such that
  \begin{equation}\label{eq:sep-seq}
    \varliminf_{n\to+\infty}\PAR{M_n(\te'_n)-M_n(\te_n)}>0.
  \end{equation}
  Then, any asymptotic M-estimators sequence $(\WH{\te}_n)_n$ is strongly
  consistent.
\end{lem}

\begin{proof}
  Suppose that \eqref{eq:sep-seq} holds for some a.s. set $A$, and that
  $(\WH{\te}_n)_n$ is a sequence of asymptotic M-estimators which is not
  strongly consistent. Saying that $(\WH{\te}_n)_n$ is not strongly consistent
  means that for any $\dP$--a.s. set, there exists a neighborhood $U$ of
  $\te^*$ and a subsequence $(\WH{\te}_{n_k})_k$ in $U^c$. In particular, on
  the a.s. set $A$, this gives a neighborhood $U$ of $\te^*$ and a
  subsequence $(\WH{\te}_{n_k})_k$ in $U^c$. Now, by virtue of
  \eqref{eq:sep-seq},
  $$
  \dP\text{--a.s},\quad \exists (\te'_{n_k})_k\in\Te^\dN,\quad
  \varliminf_{k\to+\infty}\PAR{M_{n_k}(\te'_{n_k})-M_{n_k}(\WH{\te}_{n_k})}>0,
  $$
  where the a.s. set is $A$. This contradicts \eqref{eq:ame-cns} which holds
  $\dP$--a.s. too.
\end{proof}

\begin{lem}[The $a^*$ map]\label{le:astar-prop-is-enough}
  Assume that $\Te$ is compact and that there exists a map 
  $a^*:\Te\to\Te$ such that for any
  $\te\neq\te^*$, there exists a neighborhood $U_{\te}$ of $\te$ such that
  \begin{equation}\label{eq:sep-pf}
    \dP\text{--a.s.},\ 
    \varliminf_{n\to+\infty}\inf_{U_{\te}}\PAR{M_n(a^*)-M_n}>0.
  \end{equation}
  Then, any asymptotic M-estimators sequence $(\WH{\te}_n)_n$ is strongly
  consistent.
\end{lem}

\begin{proof}
  Let us show that the assumptions of Lemma \ref{le:consistency-from-sep-seq}
  are fulfilled. We will establish \eqref{eq:sep-seq} for an a.s. set $A$
  which does not depend on the neighborhood $U$ of $\te^*$. Namely, let $U$
  be an open neighborhood of $\te^*$. For any $\te\in U^c$, let $U_\te$ and
  $A_\te$ be the neighborhood of $\te$ and the $\dP$--a.s. set for which
  \eqref{eq:sep-pf} holds. Notice that $A_\te$ depends on $U_\te$. The set
  $U^c\subset\cup_{\te\in U^c}U_{\te}$ is compact as a closed subset of the
  compact set $\Te$. We can thus extract a finite sub-covering
  $U^c\subset\cup_{i=1}^k U_{\te_i}$, and write
  \begin{align*}
    \varliminf_n \inf_{U^c}\PAR{M_n(a^*)-M_n}
    &\geq \nonumber
    \varliminf_n \min_{1\leq i\leq k}
    \inf_{U_{\te_i}}\PAR{M_n(a^*)-M_n} \\ 
    &= \nonumber
    \min_{1\leq i\leq k} \varliminf_n 
    \inf_{U_{\te_i}}\PAR{M_n(a^*)-M_n}.
  \end{align*}
  By virtue of \eqref{eq:sep-pf} we get from the above that
  \begin{equation}\label{eq:sep-U}
    \dP\text{--a.s.},\ \varliminf_n \inf_{U^c}\PAR{M_n(a^*)-M_n}>0,
  \end{equation}
  where the $\dP$--a.s. set is $A_U:=\cap_{i=1}^k A_{\te_i}$. Recall that $U$
  was a freely chosen neighborhood of $\te^*$. Consider now a countable base
  $(U_k)_k$ for $\te^*$. Then \eqref{eq:sep-U} holds on the $\dP$--a.s. set
  $A:=\cap_{i=1}^\infty A_{U_k}$, which does not depend on $U$. Notice at this
  step that
  $$
  M_n(a^*(\te_n))-M_n(\te_n) \geq \inf_{U^c}\PAR{M_n(a^*)-M_n}
  $$
  as soon as $\te_n\in U^c$ by definition of the infimum. This gives
  \eqref{eq:sep-seq} from \eqref{eq:sep-U} on the $\dP$--a.s. set $A$ defined
  above, with $(\te_n')_{n\in\dN}=(a^*(\te_n))_{n\in\dN}$.
\end{proof}

\begin{proof}[Proof of Theorem \ref{th:consistency-a-la-pf}]
  The desired result follows from Lemma \ref{le:astar-prop-is-enough}. Namely,
  let us show that \eqref{eq:sep-pf} is a consequence of \HYP{cont} and
  \HYP{sep}. Let $\te\neq\te^*$ and let $a^*$ and $V$ as in \HYP{sep}. Let
  $V_k\searrow\{\te\}$ be a decreasing local base with $V_0\subset V$. Let
  $Z:=\inf_{V}(m_{a^*}-m)$ and $Z_k:=\inf_{V_k}(m_{a^*}-m)$ and
  $Z_\infty:=m_{a^*(\te)}-m_\te$. By \HYP{cont} and the continuity of $a^*$
  and the separability of $\Te$, we get that $Z_k:\cX\to\OL{\dR}$ is
  measurable, and that
  $$
  \dP^*\text{--a.s.,} \ Z \leq Z_k \nearrow Z_\infty.
  $$
  Now, by \HYP{sep}, we get that $Z\in\rL_-^1(\cX,P^*)$ and
  $Z_\infty\in\rL_-^1(\cX,P^*)$ and $P^*(Z_\infty)>0$. Observe that
  $Z\geq-Z^-\in\rL^1(\cX,P^*)$. Thus, by the monotone convergence Theorem,
  $$
  \lim_k P^*\PAR{Z_k} = P^*(Z_\infty)>0.
  $$ 
  Therefore, $P^*(Z_k)>0$ for some $k$ (actually for $k$ large enough). Let us
  denote $U_\te:=V_k$. Now, by the law of large numbers
  $$
  \dP\text{--a.s.}, \
  \lim_n\dP_n\PAR{\inf_{U_\te}(m_{a^*}-m)}=P^*\PAR{\inf_{U_\te}(m_{a^*}-m)}>0.
  $$
  This gives finally \eqref{eq:sep-pf} since for any $n$
  $$
  \inf_{U_\te}(M_n(a^*)-M_n)=\inf_{U_\te}\dP_n(m_{a^*}-m)\geq
  \dP_n\PAR{\inf_{U_\te}(m_{a^*}-m)}.
  $$ 
\end{proof}

\textbf{Acknowledgements.} The article benefited from the comments and
criticism of the Advisory Editor and two anonymous referees. The authors would
like also to sincerely thank Professor Jon A. Wellner who has kindly answered
to their questions during his visit in Toulouse.

\providecommand{\bysame}{\leavevmode ---\ }
\providecommand{\og}{``}
\providecommand{\fg}{''}
\providecommand{\smfandname}{et}
\providecommand{\smfedsname}{\'eds.}
\providecommand{\smfedname}{\'ed.}
\providecommand{\smfmastersthesisname}{M\'emoire}
\providecommand{\smfphdthesisname}{Th\`ese}

\begin{center}
  \hrule
\end{center}

{\footnotesize
  \noindent Djalil \textsc{Chafa\"\i}, corresponding author. \\
  \textbf{Address:} UMR 181 INRA/ENVT Physiopathologie et Toxicologie
  Exp\'erimentales, \\
  \'Ecole Nationale V\'et\'erinaire de Toulouse, \\
  23 Chemin des Capelles, F-31076, Toulouse \textsc{Cedex} 3, France. \\
  \textbf{E-mail:} \url{mailto:d.chafai(AT)envt.fr} \\
  \textbf{Address:} UMR 5583 CNRS/UPS Laboratoire de Statistique et
  Probabilit\'es, \\
  Institut de Math\'ematiques de Toulouse, Universit\'e Paul Sabatier, \\
  118 route de Narbonne, F-31062, Toulouse, \textsc{Cedex} 4, France.\\
  \textbf{E-mail:} \url{mailto:chafai(AT)math.ups-tlse.fr} \\
  \textbf{Web:} \url{http://www.lsp.ups-tlse.fr/Chafai/} 
}

\end{document}